\documentclass[12pt,a4paper]{article}
\usepackage[english]{babel}
\usepackage{amssymb}
\usepackage{amsfonts,amsmath,mathtools}
\usepackage{graphics}
\usepackage{lipsum}
\usepackage[ruled]{algorithm2e}
\usepackage{float}
\usepackage{url}
\usepackage{setspace}
\usepackage{subfig}
\usepackage{xcolor}
\usepackage[export]{adjustbox}
\usepackage{tikz-cd}
\usepackage{titlesec}

\usepackage[top=4cm,bottom=3.5cm,left=3.2cm,right=3.2cm,asymmetric]{geometry}

\titleformat{\section}
{\centering\normalfont\scshape}{\thesection.}{0.5em}{}
\titleformat{\subsection}
{\normalfont}{\thesubsection.}{0.5em}{\textbf}

\usepackage{bm}
\usepackage{enumitem}
\usepackage{faktor}
\usepackage{comment}

\newtheorem{Lem}{Lemma}[section]
\newtheorem{Prop}[Lem]{Proposition}
\newtheorem{Cor}[Lem]{Corollary}
\newtheorem{Thm}[Lem]{Theorem}

\newtheorem{Def}[Lem]{Definition}
\newtheorem{Rem}[Lem]{\rm{\emph{Remark}}}

\newtheorem{Expl}[Lem]{Example}

\newenvironment{proof}[1][Proof]{\textrm{\em #1.} }{\hfill$\Box$\medskip\medskip}

\newcommand\Tor{\operatorname{Tor}}

\newcommand\Gap{\operatorname{Gap}}
\newcommand\gap{\operatorname{-gap}}

\newcommand\wdt{{\operatorname{wd}}}

\newcommand\indeg{{\operatorname{indeg}}}

\newcommand\supp{{\operatorname{supp}}}

\newcommand\lex{{\operatorname{lex}}}
\newcommand\slex{{\operatorname{slex}}}

\newcommand\Imag{{\operatorname{Im}}}

\def\L{\mathcal{L}}

\def\NZQ{\mathbb}
\def\NN{{\NZQ N}}

\def\FF{{\NZQ F}}

\let\emptyset\varnothing
\let\epsilon\varepsilon

\begin{document}
\title{\bf\normalsize\MakeUppercase{Shifting operations and completely} $t$--\MakeUppercase{spread lexsegment ideals}}
	
\author{Antonino Ficarra, Marilena Crupi}	

\newcommand{\Addresses}{{
 \footnotesize
 \textsc{Department of Mathematics and Computer Sciences, Physics and Earth Sciences, University of Messina, Viale Ferdinando Stagno d'Alcontres 31, 98166 Messina, Italy}
\begin{center}
 \textit{E-mail addresses}: \texttt{antficarra@unime.it; \texttt{mcrupi@unime.it}}.
\end{center}

}}
\date{}
\maketitle
\Addresses

\begin{abstract}
In this paper we introduce the concepts of arbitrary $t$--spread lexsegments and of arbitrary $t$--spread lexsegment ideals with $t$ a positive integer. These concepts are a natural generalization of arbitrary lexsegments and arbitrary lexsegment ideals. An ideal generated by an arbitrary $t$--spread lexsegment is called completely $t$--spread lexsegment if it is equal to the intersection of an initial $t$--spread lexsegment ideal and of a final $t$--spread lexsegment ideal. We study the class of arbitrary $t$--spread lexsegment ideals. In particular, 
we characterize all completely $t$--spread lexsegment ideals. Moreover, we classify all completely $t$--spread lexsegment ideals with a linear resolution.\\
 \emph{Keywords:} monomial ideals, minimal graded resolution, Koszul complex, $t$--spread ideals.\\
\emph{2020 Mathematics Subject Classification:} 05E40, 13B25, 13D02, 16W50.
\end{abstract}

\section{Introduction}

In the theory of graded algebras, the Hilbert series of a graded ideal is a useful tool that can detect many important properties. The Hilbert series is the generating function of the Hilbert polynomial. Let $K$ be a field.
Macaulay's Theorem \cite[Theorem 6.3.8]{JT} characterizes the numerical functions which are Hilbert functions of standard graded $K$--algebras. The proof of this result is usually achieved by studying the features of a special class of monomial ideals: the lexsegment ideals.\smallskip

Let $S = K[x_1, \ldots, x_n]$ be a polynomial ring. In the usual terminology, a monomial ideal $I$ of $S$ is a lexicographic ideal if for all monomials $u\in I$ and all monomials $v\in S$ with $\deg(u)=\deg(v)$ and such that $v\ge_{\lex}u$, then $v\in I$, where $\ge_{\lex}$ is the usual lexicographic order \cite{JT}. In this paper we call such ideals \emph{initial lexsegment ideals}. There is a rich theory that involves initial lexsegment ideals. In \cite{EK}, Eliahou and Kervaire obtained an explicit minimal free resolution for a class of ideals that includes the lexicographic one. Furthermore, Bigatti \cite{AMB} and Hulett \cite{HAH} proved, independently, that among all ideals with a given Hilbert function, the initial lexsegment ideals have the biggest Betti numbers. These results and many others have been generalized for the class of initial squarefree lexsegment ideals (see for instance \cite{JT} and the reference therein).\\
In \cite{EHQ}, a generalization of (squarefree) monomial ideals was given. Let $t\ge0$ be an integer. We say that a monomial $u=x_{i_1}x_{i_2}\cdots x_{i_d}$ is \textit{$t$--spread} if $i_{j+1}-i_j\ge t$, for all $j=1,\dots,d-1$. Every monomial is 0--spread. If $t=1$, a $t$--spread monomial is a squarefree monomial. A monomial ideal is $t$--spread if it is generated by $t$--spread monomials.\\
Let $M_{n,d,t}$ be the set of all $t$--spread monomials of $S$ of degree $d$, $t\ge 1$ and assume that $M_{n,d,t}$ is endowed with the squarefree lexicographic order $\ge_{\slex}$. Let $u,v\in M_{n,d,t}$. If $u\ge_{\slex}v$, the set $\mathcal{L}_t(u,v)=\{\omega\in M_{n,d,t}:u\ge_{\slex}\omega\ge_{\slex}v\}$ is called an \emph{arbitrary $t$--spread lexsegment}.
If $u$ is the maximum monomial of $M_{n,d,t}$, or $v$ is the minimum monomial of $M_{n,d,t}$, we say that $\mathcal{L}_t(u,v)$ is an \emph{initial} or a \emph{final $t$--spread lexsegment}, respectively. The ideal generated by an arbitrary $t$--spread lexsegment is called an \emph{arbitrary $t$--spread lexsegment ideal}.
Let $I = (\mathcal{L}_t(u,v))$ be an arbitrary $t$--spread lexsegment ideal, $I$ is said \emph{completely $t$--spread lexsegment} if $I=(\mathcal{L}_t(\max(M_{n,d,t}),v))\cap(\mathcal{L}_t(u,\min(M_{n,d,t})))$. This definition agrees with the one given by De Negri and Herzog (\cite[Proposition 1.4 (a)]{ADH}) in all cases but the final $t$--spread lexsegment case.\\
Every initial and every final $t$--spread lexsegment ideal is completely $t$--spread lexsegment, but not all arbitrary $t$--spread lexsegment ideals are completely $t$--spread lexsegment ideals. In this paper we characterize this latter class. The case $t=0$ has been deeply treated by De Negri and Herzog in \cite{DH}. In this paper we consider the case $t\ge1$. To characterize the completely $t$--spread lexsegment ideals we use the following useful decomposition:\medskip
\\ (*) \, given $t\ge1$, $u,v\in M_{n,d,t}$ with $u\ge_{\slex}v$, then, setting $N=\mathcal{L}_t\big(\max(M_{n,d,t}), v\big)$, $F=\mathcal{L}_t\big(u,\min(M_{n,d,t})\big)$ and $L=\mathcal{L}_t(u,v)$,
one has $L=N\cap F$.\medskip \\
The aim of this article is to characterize all completely $t$--spread lexsegment ideals for $t\ge 1$. Some shifting operations will be a fundamental tool for the development of the paper. \\
The plan of the paper is as follows. Section \ref{sec1} contains some preliminary materials. In Section \ref{sec2}, we prove our main result which is a characterization of the class of completely $t$--spread lexsegment ideals (Theorem \ref{CharacterizationCompletelyLexSegTSpread}). For its proof we need some preliminary lemmas which describe a basis of the Koszul homologies of some suitable classes of arbitrary $t$--spread ideals (Lemma \ref{lemmakoszulhomtspread}). We point out that a final $t$--spread lexsegment is a $t$--spread strongly stable ideal but with the order on the variables reversed. In Section \ref{sec3} we determine all completely $t$--spread lexsegment ideals with a linear resolution. The paper contains many examples illustrating the main results. 

\section{Preliminaries}\label{sec1}

Let $K$ be a field and let $S=K[x_1,\dots,x_n]$ be the standard graded polynomial ring in $n$ indeterminates. $S$ is a $\NN$--graded $K$--algebra with each variable of degree $1$. For any monomial ideal $I$ of $S$, there exists a unique minimal graded free $S$--resolution of the form
$$
\FF: 0\rightarrow \bigoplus_{j\ge0}S(-j)^{\beta_{r,j}}\xrightarrow{\ d_r\ } \ldots\xrightarrow{\ d_1\ }\bigoplus_{j\ge0}S(-j)^{\beta_{0,j}}\xrightarrow{\ d_0\ }I\rightarrow0,
$$
where $S(-j)$ denotes the free $S$--module obtained by shifting the degrees of $S$ by $j$ and $r\le n$. The numbers $\beta_{i,j}=\beta_{i,j}(I) = \dim_K\Tor^S_i(K,I)_j$ are called the \textit{graded Betti numbers} of $I$.
We say that $I$ has a \textit{$d$--linear resolution} if $\beta_{i,j}(I) =0$, for all $i>0$ and $j\neq d+i$. It follows that $I$ is generated in degree $d$. Equivalently, $I$ has a \textit{$d$--linear resolution} if $\Tor^S_i(K,S/I)_j=0$, for all $i>0$ and all $j\ne d+i-1$.

If $I$ is a monomial ideal of $S$, we denote by  $G(I)$ the unique minimal set of monomial generators of $I$, whereas $\indeg(I)$ is the minimum degree of the monomials of $G(I)$. For a monomial $u\in S$, $u\ne 1$, we set \[\supp(u)=\big\{i:x_i\ \text{divides}\ u\big\},\]
and write
\[\max(u)=\max\big\{i:i\in\supp(u)\big\},\,\,
\min(u)=\min\big\{i:i\in\supp(u)\big\}.\]
Moreover, we set $\max(1)=\min(1)=0$.

In \cite{EHQ}, the following definition has been introduced.
\begin{Def}
	\rm Given two integers $n\ge1,t\ge0$, a monomial $x_{i_1}x_{i_2}\cdots x_{i_d}\in S= K[x_1, \ldots, x_n]$, with $1\le i_1\le i_2\le\dots\le i_d\le n$, is said \textit{$t$--spread} if $i_{j+1}-i_j\ge t$, for all $j=1,\dots,d-1$.\\
	We say that a monomial ideal $I$ of $S$ is \textit{$t$--spread} if it is generated by $t$--spread monomials.
\end{Def}

For instance, $x_2x_4x_6\in K[x_1,x_2,x_3,x_4,x_6]$ is a $2$--spread monomial, but not a $3$--spread monomial.  One can note that any monomial is $0$--spread and any monomial ideal is a $0$--spread monomial ideal. Moreover, any squarefree monomial is $1$--spread and any squarefree monomial ideal is a $1$--spread ideal.\\

For the reader's convenience, we quote some notations and results from \cite{CAC}. 

Let $n,d,t$ be positive integers with $n\ge 1+(d-1)t$ and denote by $M_{n,d,t}$ the set of all $t$--spread monomials of $S$ of degree $d$. From \cite[Theorem 2.3]{EHQ}, we have that 
\[|M_{n,d,t}|=\binom{n-(d-1)(t-1)}{d}.\] 

The monomial ideal generated by $M_{n,d,t}$ is denoted by $I_{n,d,t}$ \cite{EHQ} and is called the \textit{$t$--spread Veronese ideal of degree $d$} of the polynomial ring $S=K[x_1,\dots,x_n]$.\\

Let $t\ge 1$. Throughout this paper $M_{n,d,t}$ is endowed with the \emph{squarefree lexicographic order}, $\ge_{\slex}$, \cite{AHH2}. Given $u,v\in M_{n,d,t}$ with $u=x_{i_1}x_{i_2}\cdots x_{i_d}$, $1\le i_1<i_2<\dots<i_d\le n$, $v=x_{j_1}x_{j_2}\cdots x_{j_d}$, 
$1\le j_1<j_2<\dots<j_d\le n$, 
then $u>_{\slex}v$ if $i_1=j_1,\dots,i_{s-1}=j_{s-1}$ and $i_s<j_s$ for some $1\le s\le d$. 

Let $U$ be a not empty subset of $M_{n,d,t}$. We denote by $\max(U)$ ($\min(U)$, respectively) the maximum (minimum, respectively) monomial $w\in U$, with respect to $>_{\slex}$.

\begin{Def}\rm
Let $u,v\in M_{n,d,t}$, $u\ge_{\slex}v$.
The set
\[\mathcal{L}_t(u,v) = \big\{\omega\in M_{n,d,t}:u\ge_{\slex}\omega\ge_{\slex}v\big\}\]
is called an \emph{arbitrary $t$--spread lexsegment set}, or simply a \emph{$t$--spread lexsegment set}. The set
\[\mathcal{L}_t^{i}(v) = \big\{\omega\in M_{n,d,t}:\omega\ge_{\slex}v\big\} = \mathcal{L}_t(\max(M_{n,d,t}),v) = \mathcal{L}_t(x_1x_{1+t}\cdots x_{1+(d-1)t},v)\]
is called an \emph{initial $t$--spread lexsegment} and 
\[\mathcal{L}_t^{f}(u)=\big\{\omega\in M_{n,d,t}:\omega\le_{\slex}u\big\}= \mathcal{L}_t(u,\min(M_{n,d,t}))= \mathcal{L}_t(u,x_{n-(d-1)t}\cdots x_{n-t}x_n)\]
is called a \emph{final $t$--spread lexsegment}.
\end{Def}

One can observe that, if $u,v\in M_{n,d,t}$, with $u\ge_{\slex}v$, then $\mathcal{L}_t(u,v)=\mathcal{L}_t^i(v)\cap\mathcal{L}_t^f(u)$.  This observation will be pivotal in the sequel. 

Let $n,d,t\ge 1$ be positive integers with $n\ge 1+(d-1)t$. 
\begin{Def}\label{def:comp}\rm
	A $t$ spread graded ideal $I$ of $S$ is called an \textit{arbitrary $t$--spread lexsegment ideal}, or simply a \textit{$t$--spread lexsegment ideal} if it is generated by an (arbitrary) $t$--spread lexsegment. 
	\\
	A $t$--spread lexsegment ideal $I=(\mathcal{L}_t(u,v))$ is a \textit{completely $t$--spread lexsegment ideal} if
	$$
	I=J\cap T,
	$$
	where $J=(\mathcal{L}_t^i(v))$ and $T=(\mathcal{L}_t^f(u))$.
\end{Def}

More in detail, a $t$--spread ideal $I$ of $S$ is a $t$--spread lexsegment ideal if it is generated in one degree $d$ and $I=(\mathcal{L}_t(u,v))$, with $u,v\in M_{n,d,t}$, $u\ge_{\slex}v$. Furthermore, the ideal $I$ is a \textit{completely $t$--spread lexsegment ideal} if it is generated in one degree $d$ and furthermore $I=(\mathcal{L}_t^i(v))\cap(\mathcal{L}_t^f(u))$.

\begin{Expl} \label{ex1}\rm Let $S= K[x_1, \ldots, x_{11}]$. Consider the $3$--spread lexsegment ideal $I = (\mathcal{L}_t(u,v))$, with  $u=x_1x_5x_8$ and $v=x_2x_5x_8$. We have:
	\begin{align*}
	I=&(x_{1}x_{5}x_{8},x_{1}x_{5}x_{9},x_{1}x_{5}x_{10},x_{1}x_{5}x_{11},x_{1}x_{6}x_{9},
	x_{1}x_{6}x_{10},\\& \ x_{1}x_{6}x_{11}, x_{1}x_{7}x_{10},x_{1}x_{7}x_{11},x_{1}x_{8}x_{11},x_{2}x_{5}x_{8}),\\
	J=&  (\mathcal{L}_t^i(v))=(x_{1}x_{4}x_{7},x_{1}x_{4}x_{8},x_{1}x_{4}x_{9},x_{1}x_{4}x_{10},x_{1}x_{4}x_{11},x_{1}x_{5}x_{8},\\
	&\ x_{1}x_{5}x_{9},x_{1}x_{5}x_{10},x_{1}x_{5}x_{11},x_{1}x_{6}x_{9},x_{1}x_{6}x_{10},x_{1}x_{6}x_{11},\\
	&\ x_{1}x_{7}x_{10},x_{1}x_{7}x_{11},x_{1}x_{8}x_{11},x_{2}x_{5}x_{8}),\\
	T=&(\mathcal{L}_t^f(u)) = (x_{1}x_{5}x_{8},x_{1}x_{5}x_{9},x_{1}x_{5}x_{10},x_{1}x_{5}x_{11},x_{1}x_{6}x_{9},x_{1}x_{6}x_{10}, x_{1}x_{6}x_{11},\\
	&\ x_{1}x_{7}x_{10},x_{1}x_{7}x_{11},x_{1}x_{8}x_{11},x_{2}x_{5}x_{8},x_{2}x_{5}x_{9},x_{2}x_{5}x_{10},x_{2}x_{5}x_{11},\\
	&\ x_{2}x_{6}x_{9},x_{2}x_{6}x_{10},x_{2}x_{6}x_{11},x_{2}x_{7}x_{10},x_{2}x_{7}x_{11}, x_{2}x_{8}x_{11},x_{3}x_{6}x_{9},\\ &\ x_{3}x_{6}x_{10},x_{3}x_{6}x_{11},x_{3}x_{7}x_{10},x_{3}x_{7}x_{11}, x_{3}x_{8}x_{11},x_{4}x_{7}x_{10},x_{4}x_{7}x_{11},\\ &\ x_{4}x_{8}x_{11},x_{5}x_{8}x_{11}).
	\end{align*}
One can quickly verify that $I=J\cap T$, \emph{i.e.},  $I$ is completely $3$--spread lexsegment.
\end{Expl}
\begin{Expl} \label{ex:nocompl}\rm Let $I = (x_1x_5x_7, x_2x_4x_6)$ be a 2--spread ideal of $S=K[x_1, \ldots, x_7]$. $I$ is a $2$--spread lexsegment ideal but not a completely $2$--spread lexsegment ideal. Indeed, setting $u = x_1x_5x_7$ and $v = x_2x_4x_6$, then 
\[J = (\mathcal{L}_t^{i}(v)) = (x_1x_3x_5, x_1x_3x_6,x_1x_3x_7, x_1x_4x_6, x_1x_4x_7, x_1x_5x_7, x_2x_4x_6)\]
\[T = (\mathcal{L}_t^{f}(u)) = (x_1x_5x_7, x_2x_4x_6, x_2x_4x_7, x_2x_5x_7, x_3x_5x_7),\]
and $J \cap T = (x_1x_5x_7, x_2x_4x_6, x_1x_2x_4x_7)\neq I$. Note that $J \cap T$ is not a $2$--spread ideal of $S$.
\end{Expl}
\begin{Rem} \em Definition \ref{def:comp} is different from the one given in \cite{DH} by De Negri and Herzog. More precisely, their definition involves the notion of \emph{shadow} and is equivalent to our definition in some cases (\cite[Theorem 2.3]{DH}, \cite[Proposition 1.5]{ADH}). Our definition is forced by the fact that the intersection of $t$--spread ideals does not give always a $t$--spread ideal (see, for instance Example \ref{ex:nocompl}).
\end{Rem}

In this paper we want to characterize all completely $t$--spread lexsegment ideals. 

One can observe that the $t$--spread Veronese ideal $I_{n,d,t}$ is the unique $t$--spread lexsegment ideal of degree $d$ of the polynomial ring $S$ that is both an initial and a final $t$--spread lexsegment ideal, thus a completely $t$--spread lexsegment ideal. Moreover, an initial (final) $t$--spread lexsegment ideal $I$ of $S$ is a completely $t$--spread lexsegment ideal too. Indeed, $I = I\cap I_{n,d,t}$ in both cases.\\

For later use, we recall some notions and properties from \cite{EHQ}.

A subset $U$ of $M_{n,d,t}$ is called a \emph{$t$--spread strongly stable set}, if the following property holds: for all $u\in U$, all $j\in \supp(u)$ and all $1\le i< j$ such that $x_i(u/x_j)\in M_{n,d,t}$, it follows that $x_i(u/x_j)\in U$.

\begin{Def}\rm
	A $t$--spread monomial ideal $I$ is a \textit{$t$--spread strongly stable ideal} if for all $t$--spread monomial $u\in I$, all $j\in \supp(u)$ and all $1\le i< j$ such that $x_i(u/x_j)\in M_{n,d,t}$, it follows that $x_i(u/x_j)\in I$.
\end{Def}

Every initial $t$--spread lexsegment ideal $I$ of $S$ is a $t$--spread strongly stable ideal. Hence, in order to compute the graded Betti numbers of $I$ one may use the Ene, Herzog, Qureshi's formula (\cite[Corollary 1.12]{EHQ}):
\begin{equation}
\label{eq1}
\beta_{i,i+j}(I)\ =\ \sum_{u\in G(I)_j}\binom{\max(u)-t(j-1)-1}{i}.
\end{equation}
For $t=0$, the previous formula is the Eliahou--Kervaire formula \cite{EK} for (strongly) stable ideals, whereas for $t=1$, it gives
the Aramova--Herzog--Hibi formula \cite{AHH2, JT}.\\

A $t$--spread final lexsegment ideal $I$ of $S$ is also a $t$--spread strongly stable ideal, but with the order of the variables reversed, \emph{i.e.}, $x_n> x_{n-1} > \cdots >x_1$. Hence, in order to compute their graded Betti numbers one may use the following \emph{modified} Ene, Herzog, Qureshi's formula:
\begin{equation}\label{Bettinumbersreversetspread}
\beta_{i,i+j}(I)\ =\ \sum_{u\in G(I)_j}\binom{n-\min(u)-t(j-1)}{i}.
\end{equation}

One can note that both an initial and a final $t$--spread final lexsegment ideal have a linear resolution since they are $t$--spread strongly stable ideal generated in one degree (\cite[Theorem 1.4]{EHQ}). On the other hand, one can obtain such a property directly from (\ref{eq1}) and (\ref{Bettinumbersreversetspread}), respectively.

We close the section by recalling some details on the \textit{Koszul complex}. Let ${\bf x}=x_1,\dots,x_n$ a sequence of elements of $S$. The \textit{Koszul complex} $K_{_{\text{\large$\boldsymbol{\cdot}$}}}({\bf x};S)$ attached to the sequence ${\bf x}$ is defined as follows: let $F$ be the free $S$--module with basis $e_1,\dots,e_n$. 
\begin{enumerate}
	\item[--] We let $K_j({\bf x};S)={\bigwedge}^jF$, for all $j=0,\dots,n$. A basis of the free $S$--module $K_j({\bf x};S)$ is given by the wedge products $e_\sigma=e_{i_1}\wedge e_{i_2}\wedge\dots\wedge e_{i_j}$, where $\sigma=\{i_1<i_2<\dots<i_j\}\subseteq\{1,\dots,n\}$.
	\item[--] We define the differential $\partial_j:K_j({\bf x};S)\rightarrow K_{j-1}({\bf x};S)$, $j=1,\dots,n-1$ by 
	$$
	\partial_j(e_\sigma)=\sum_{k=1}^j(-1)^{k+1}x_{i_k}e_{\sigma\setminus\{i_k\}}.
	$$
\end{enumerate}
If $I$ is a graded ideal of $S$, one can define the complex $K_{_{\text{\large$\boldsymbol{\cdot}$}}}({\bf x};S/I) = K_{_{\text{\large$\boldsymbol{\cdot}$}}}({\bf x};S)\otimes S/I$. We denote by $H_i({\bf x};S/I)= H_i(K_{_{\text{\large$\boldsymbol{\cdot}$}}}({\bf x};S/I))$ the $i$--th homology module of ${\bf x}$ with respect to $S/I$.\\
Since, ${\bf x}$ is a regular sequence on $S$, there is an isomorphism of graded $K\cong S/({\bf x})$--vector spaces \cite[Corollary A.3.5]{JT}
$$
\Tor^S_i(K,S/I)\cong H_i({\bf x};S/I),
$$
and so $\beta_{i,j}(S/I)=\dim_K H_i({\bf x};S/I)_j$, for all $i,j\ge0$.\\

\section{A characterization of completely $t$--spread lexsegment ideals}\label{sec2}

In this section we give a characterization of a completely $t$--spread lexsegment ideal $I$ of the polynomial ring $S = K[x_1, \ldots, x_n]$. From now on, given three positive integers $n,d,t$, we tacitly always assume that $n\ge 1+(d-1)t$, otherwise $M_{n,d,t}=\emptyset$.

We start recalling a notion introduced  in \cite[Definition 3.7]{AFC2} (see also \cite{AC2, AC1}).

Given an integer $m$, by $[m]$ we denote the set $\{1,\dots,m\}$. Let $u=x_{i_1}x_{i_2}\cdots x_{i_d}$ be a $t$--spread monomial of degree $d$ of $S$, with $1\le i_1<\dots<i_d\le n$. Let $j\in[d-1]$. We say that $u$ has a \textit{$j\gap_t$} if $i_{j+1}-i_j-t>0$. For all $j=1,\dots,d-1$, we set 
$$
\wdt(j\gap_t)(u)=i_{j+1}-i_j-t
$$
and we call it the \textit{width of the $j\gap_t$}.

Finally, we denote by $\Gap_t(u)$ the set defined as follows:
$$
\Gap_t(u) = \big\{j\in[d-1]:\text{there exists a}\ j\gap_t(u)\big\}.
$$
We can now reformulate for our purpose a result given in \cite{AFC1} \emph{via} gaps.
\begin{Lem}
	\label{maxbettilemma}
	\textup{\cite[Lemma 3.1]{AFC1}} Given two integers $n, t\ge1$, let $u=x_{i_1}x_{i_2}\cdots x_{i_d}$, with $1\le i_1<i_2<\dots<i_d\le n$, be a $t$--spread monomial of degree $d$ of $S=K[x_1,\dots,x_n]$. If $\Gap_t(u)=\emptyset$, then $u$ is the smallest $t$--spread monomial of degree $d$ of $S$, with respect to $>_{\slex}$. 
	Otherwise, if $p=\max(\Gap_t(u))$, then the largest $t$--spread monomial $v$ of degree $d$ of $S$ following $u$, with respect to $>_{\slex}$, is
	$$
	v = x_{i_1}\cdots x_{i_{p-1}}x_{i_p+1}x_{i_p+1+t}\cdots x_{i_p+1+t(d-p-1)}x_{i_p+1+t(d-p)}.
	$$
\end{Lem}

\begin{Rem} \rm Lemma \ref{maxbettilemma} is very useful. Indeed, given $u,v\in M_{n,d,t}$, $t\ge1$ such that $u\ge_{\slex}v$, it allows to determine all the $t$--spread monomials in the set $L=\mathcal{L}_t(u,v)$.
\end{Rem}

Now, let $n,d,t\ge1$. In \cite{EHQ}, the following two operators 
\begin{align*}
\sigma^{t-1}:x_{i_1}x_{i_2}\cdots x_{i_d}\in M_{n,d,1}&\longmapsto\prod_{j=1}^dx_{i_j+(j-1)(t-1)}\in M_{n+(d-1)(t-1),d,t},\\
\tau^{t-1}:x_{i_1}x_{i_2}\cdots x_{i_d}\in M_{n,d,t}&\longmapsto\prod_{j=1}^dx_{i_j-(j-1)(t-1)}\in M_{n-(d-1)(t-1),d,1},
\end{align*}
have been introduced.
The operators $\sigma^{t-1}$ and $\tau^{t-1}$ are both bijections. 

Applying Lemma \ref{maxbettilemma} we obtain the next result.

\begin{Cor}\label{cortaut-1lexseg}
	Let $n,d,t\ge1$ be positive integers and let $u,v$ be $t$--spread monomials of $S$ with $u\ge_{\slex}v$. Setting $L=\mathcal{L}_t(u,v)$, then $$\tau^{t-1}(L)=\mathcal{L}_1(u^{\tau^{t-1}},v^{\tau^{t-1}})$$ is a squarefree lexsegment in the polynomial ring $K[x_1,\dots,x_{n-(d-1)(t-1)}]$.
\end{Cor}
\begin{proof}
	Firstly, one can note that if $y=\min(M_{n,d,t}) = x_{n-(d-1)t}\cdots x_{n-t}x_n$,
	then
	\[y^{\tau^{t-1}}=x_{n-(d-1)t}x_{n-(d-1)t+1}\cdots x_{n-(d-1)(t-1)}= \min(M_{n-(d-1)(t-1),d,1}).\]
	Now, let $y\in M_{n,d,t}$, $y\neq \min(M_{n,d,t})$ and let $\omega=\max\{z\in M_{n,d,t}:z<_{\slex}y\}$, then
	$$
	\omega^{\tau^{t-1}}=\max\big\{z\in M_{n-(d-1)(t-1),d,1}:z<_{\slex}y^{\tau^{t-1}}\big\}.
	$$
	In fact,  one can quickly verify that $\max(\Gap_t(y))=\max(\Gap_1(y^{\tau^{t-1}}))$. The assertion follows from Lemma \ref{maxbettilemma}.
\end{proof}

\begin{Prop}\label{propgendandd+1}
	Let $n,d$ be positive integers, $t\ge0$. Let $u,v\in M_{n,d,t}$, $u\ge_{\slex}v$, $J=(\mathcal{L}_t^i(v))$ and $T=(\mathcal{L}_t^f(u))$. Then $Q=J\cap T$ is generated in degrees $d$ and $d+1$.
\end{Prop}
\begin{proof}
Since $J$ is a $t$--spread lexsegment initial ideal generated in degree $d$, then it 
has a $d$--linear resolution, and so
\begin{equation}\label{torIlinearres}
\Tor_{i}^S(K,S/J)_j\ =\ 0,\ \ \ \text{for all}\ \ i>0\ \ \text{and all}\ \ j\ne i+d-1.
\end{equation}
Similarly, since $T$ has a $d$--linear resolution generated in degree $d$, one has
\begin{equation}\label{torJlinearres}
\Tor_{i}^S(K,S/T)_j\ =\ 0,\ \ \ \text{for all}\ \ i>0\ \ \text{and all}\ \ j\ne i+d-1.
\end{equation}
Now, set $V=J+T$. It is clear that $V=I_{n,d,t}$ is the $t$--spread Veronese ideal of degree $d$ of $S$. The
short exact sequence of $S$--modules
$$
0\rightarrow S/Q\rightarrow S/J\oplus S/T\rightarrow S/V\rightarrow 0
$$
induces the long exact sequence of $\Tor$--modules:
\begin{equation}\label{longexactsequencetortspread}
\begin{aligned}
	\phantom{\cdots \rightarrow \hfill\Tor_{i}^S(K,S/Q)}& \phantom{\rightarrow \Tor_{i}^S(K,S/J..*.}\cdots\rightarrow\hfill\Tor_{i+1}^S(K,S/V) \rightarrow \hphantom{0}\\
	\rightarrow \hfill\Tor_{i}^S(K,S/Q) & \rightarrow \Tor_{i}^S(K,S/J\oplus S/T) \rightarrow\hfill\Tor_{i}^S(K,S/V) \rightarrow\hphantom{0}\\
	\rightarrow
	\hfill\Tor_{i-1}^S(K,S/Q)
	& \rightarrow \makebox[\widthof{$\Tor_{i}^S(K,S/J\oplus S/T)$}][c]{$\cdots\hfill \cdots$} \rightarrow
	\hfill\Tor_{0}^S(K,S/V) \rightarrow0.
	\end{aligned}
\end{equation}

Since $V=I_{n,d,t}$ has a $d$--linear resolution \cite[Theorem 2.3 (c)]{EHQ}, we have $\Tor_{i+1}^S(K,S/V)_j$ $=$ $0$, for all $i>0$ and all $j\ne i+d$. From this fact and by (\ref{torIlinearres}), (\ref{torJlinearres}) and (\ref{longexactsequencetortspread}) we deduce that $\Tor_{i}^S(K,S/Q)_j=0$ for all $i>0$ and all $j\ne i+d-1,i+d$. In particular, setting $\mathfrak{m}=(x_1,\dots,x_n)$, we have that 
$$
0= \Tor_{1}^S(K,S/Q)_j \cong \Big((\mathfrak{m}\cap Q)\big/\mathfrak{m}Q\Big)_j=Q_j\big/\mathfrak{m}Q_{j-1},\quad j\ne d,d+1.
$$ 
Hence $Q$ is generated in degrees $d$ and $d+1$.
\end{proof}

For the proof of the characterization, we need some technical lemmas that describe the basis of some homologies of the Koszul complex, seen as $K$--vector spaces. Following \cite{ADH}, it is convenient to denote the image of a monomial $u\in S$ in any quotient ring of $S$ again by $u$. We will keep this convention throughout this paper.

We quote next result from \cite{AHH2}.

\begin{Lem}\label{AraHerHibiKoszulHomology}
	\textup{\cite[Proposition 2.2]{AHH2}} Let $I$ be a squarefree strongly stable ideal of $S$. Then, for every $i>0$, a basis of the homology classes of $H_i(x_1,\dots,x_n;S/I)$ is given by the homology classes of the cycles
	\begin{align*}
	& u/x_{\max(u)}e_{\sigma}\wedge e_{\max(u)},
	\end{align*}
	with $u\in G(I)$, $|\sigma|=i-1$, $1\le\max(\sigma)<\max(u)$, $\sigma\cap\supp(u)=\emptyset$.
\end{Lem}

This result will be used for computing a basis of the $K$--vector space $$H_i(x_1,\dots,x_n;S/I)$$ for $I$ initial squarefree lexsegment ideal of $S$. Indeed, $I$ is a squarefree strongly stable ideal \cite{AHH2}. Following \cite{DH} (see also \cite{AH}), the cycles in Lemma \ref{AraHerHibiKoszulHomology} will be called \emph{classical cycles}.

Now we prove the $t$--spread version of the previous result.
\begin{Lem}\label{lemmakoszulhomtspread}
	Let $I$ be a $t$--spread strongly stable ideal of $S$, $t\ge 1$. For all $i>0$, a basis of the homology classes of $H_i(x_1,\dots,x_n;S/I)$ is given by the homology classes of the cycles.
	$$
	u^{\tau^{t-1}}/x_{\max(u^{\tau^{t-1}})}e_{\sigma}\wedge e_{\max(u^{\tau^{t-1}})},
	$$
	with $u\in G(I)$, $|\sigma|=i-1$, $\sigma\cap\supp(u^{\tau^{t-1}})=\emptyset$, $1\le \max(\sigma)<\max(u^{\tau^{t-1}})$.
\end{Lem}
\begin{proof}
	Let $\widetilde S=K[x_1,\dots,x_{n-(\indeg(I)-1)(t-1)}]$. Since $I$ is a $t$--spread strongly stable ideal of $S$, $I^{\tau^{t-1}}$ is a squarefree strongly stable ideal of the polynomial ring $\widetilde S$, and $\beta_{i,j}=\beta_{i,j}(I)=\beta_{i,j}(I^{\tau^{t-1}})$ for all $i,j\ge0$ (see \cite{EHQ}). Therefore,
	for all $i,j\ge0$ we have a graded isomorphism of $K$--vector spaces
	$$
	\Tor^S_i(K,S/I)_j\cong\Tor^{\widetilde{S}}_i(K,\widetilde S/I^{\tau^{t-1}})_j.
	$$
	Hence, we have an isomorphism of graded $K$--vector spaces
	\begin{equation}\label{isomTorij}
	\Tor^S_i(K,S/I)\cong \Tor^{\widetilde S}_i(K,\widetilde S/I^{\tau^{t-1}}).
	\end{equation}
	We observe that ${\bf x}=x_1,\dots,x_n$ is a regular sequence on $S$, and $({\bf x})=\mathfrak{m}$ is the maximal graded ideal of $S$. Analogously, $\widetilde{\bf x}=x_1,\dots,x_{n-(\indeg(I)-1)(t-1)}$ is a regular sequence on $\widetilde S$, and $\widetilde{\mathfrak{m}}=(\widetilde{\bf x})$ is the graded maximal ideal of $\widetilde S$. Therefore, by \cite[Theorem A.3.5]{JT}, we have the isomorphism of graded $K$--vector spaces, for all $i\ge0$,
	\begin{equation}\label{isomHiTori}
	\begin{aligned}
	H_i({\bf x};S/I)&\cong\Tor_{i}^S(K,S/I),\\
	H_i(\widetilde{\bf x}; \widetilde S/I^{\tau^{t-1}})&\cong\Tor_{i}^{\widetilde S}(K,\widetilde S/I^{\tau^{t-1}}).
	\end{aligned}
	\end{equation}
	
	Finally, by (\ref{isomTorij}) and (\ref{isomHiTori}), we have the following isomorphism of graded $K$--vector spaces
	\begin{align*}
	H_i({\bf x};S/I)\cong H_i(\widetilde{\bf x};\widetilde S/I^{\tau^{t-1}}).
	\end{align*}
	
	Hence, to compute a basis of $H_i({\bf x};S/I)$ it suffices to determine a basis of $H_i(\widetilde{\bf x};\widetilde S/I^{\tau^{t-1}})$. 
	
	Since $G(I^{\tau^{t-1}})=\big\{u^{\tau^{t-1}}:u\in G(I)\big\}$, by Lemma \ref{AraHerHibiKoszulHomology}, the basis is given by the homology classes of the Koszul cycles 
	\[
	u^{\tau^{t-1}}/x_{\max(u^{\tau^{t-1}})}e_{\sigma}\wedge e_{\max(u^{\tau^{t-1}})} 
	\]
	with $u\in G(I)$, $|\sigma|=i-1$, $\sigma\cap\supp(u^{\tau^{t-1}})=\emptyset$, $1\le\max(\sigma)<\max(u^{\tau^{t-1}})$.
	Note that $\max(u^{\tau^{t-1}})= \max(u)-(\deg(u)-1)(t-1)$.
\end{proof}

One can observe that for $t=1$, Lemma \ref{lemmakoszulhomtspread} returns Lemma \ref{AraHerHibiKoszulHomology}. Indeed, $u^{\tau^{t-1}}=u$, in such a case.
\\
Lemma \ref{lemmakoszulhomtspread} will be used for computing a basis of the $S$--module $H_i(x_1,\dots,x_n;S/I)$ for an initial $t$--spread lexsegment ideal $I$ of $S$, $t\ge 1$. In fact, such an ideal is $t$--spread strongly stable.  

As we have previously pointed out, a final $t$--spread lexsegment ideal $I$ of $S$ is a $t$--spread strongly stable ideal of $S$ with the order on the variables reversed $x_n>x_{n-1}>\dots >x_1$. Hence, a basis of the $S$--module $H_i(x_1,\dots,x_n; S/I)$ is given by the homology classes of the cycles 
\begin{equation}\label{eq2:fintspread}
	u^{\tau^{t-1}}/x_{\min(u^{\tau^{t-1}})}e_{\min(u^{\tau^{t-1}})}\wedge e_{\sigma} = u^{\tau^{t-1}}/x_{\min(u)}e_{\min(u)}\wedge e_{\sigma},
	\end{equation}
	with $u\in G(I)$, $|\sigma|=i-1$, $\sigma\cap\supp(u^{\tau^{t-1}})=\emptyset$, $1\le\min(u)<\min(\sigma)\le\max(\sigma)\le n-(\deg(u)-1)(t-1)$.

In particular, if $I$ is a final squarefree lexsegment ideal of $S$, a basis of the homology classes of $H_i(x_1,\dots,x_n;S/I)$ of $I$ is given by the homology classes of the cycles
\begin{align}\label{eq:revhom}
	&u/x_{\min(u)}e_{\min(u)}\wedge e_{\sigma},
	\end{align}
	with $u\in G(I)$, $|\sigma|=i-1$, $\min(u)<\min(\sigma)\le\max(\sigma)\le n$, $\sigma\cap\supp(u)=\emptyset$.

Following \cite{DH}, the cycles in (\ref{eq2:fintspread}) will be called \emph{non classical cycles}.

\begin{Thm}\label{CharacterizationCompletelyLexSegTSpread}
Given $n,d,t$ positive integers, let $u=x_{i_1}x_{i_2}\cdots x_{i_d},v=x_{j_1}x_{j_2}\cdots x_{j_d}$ be two $t$--spread monomials of degree $d$ of $S=K[x_1,\dots,x_n]$ such that $u\ge_{\slex}v$. Let $I$ be the ideal generated by $\L_t(u,v)$. Set
$$\quad J=(\mathcal{L}_t^i(v)), \quad T=(\mathcal{L}_t^f(u)).$$
Then, the following conditions are equivalent:
\begin{enumerate}
\item[\textup{(a)}] $I$ is a completely $t$--spread lexsegment ideal, \emph{i.e.}, $I=J\cap T$;
\item[\textup{(b)}] for every $\omega\in M_{n,d,t}$ with $\omega<_{\slex}v$ there exists an integer $s>i_1$ such that $x_s$ divides $\omega$ and $x_{i_1}(\omega/x_{s})\le_{\lex}u$. 
\end{enumerate}
\end{Thm}
\begin{proof}
Setting $I=(\L_t(u,v))$, $N = \mathcal{L}_t^i(v)$ and $F = \mathcal{L}_t^f(u)$, let us define 
\[J = (N),\quad T = (F).\]
It is $L= \mathcal{L}_t(u,v) = N\cap F$. From Proposition \ref{propgendandd+1}, $J\cap T$ is generated in degrees $d$ and $d+1$. On the other hand, $I\subseteq J\cap T$ and $L=[I_d]_t=[(J\cap T)_d]_t=N\cap F$. 
Therefore, $I$ is a completely $t$--spread lexsegment ideal, \emph{i.e.}, $I=J\cap T$, if and only if $J\cap T$ is generated in degree $d$. Set $V=J+T=I_{n,d,t}$. 

Let us consider the exact sequence (see Proposition \ref{propgendandd+1}, proof) 
$$
0\rightarrow S/J\cap T\rightarrow S/J\oplus S/T\longrightarrow S/V\rightarrow0
$$
which induces the following long exact sequence
\begin{equation}\label{exactsequencecomletelytspreadlexseg}
\begin{aligned}
& \cdots \rightarrow \hfill\Tor_{2}^S(K,S/J\cap T)_{d+1} \rightarrow \hfill
\Tor_{2}^S(K,S/J\oplus S/T)_{d+1} \xrightarrow{\ \alpha\ }\\
&\xrightarrow{\ \alpha\ }\Tor_{2}^S(K,S/V)_{d+1} \xrightarrow{\ \beta\ } \Tor_{1}^S(K,S/J\cap T)_{d+1}\rightarrow 0.
\end{aligned} 
\end{equation}
\textsc{Claim 1.}  The condition $I=J\cap T$ is equivalent to the fact that $\alpha$ is surjective.\\
Indeed, suppose that  $\alpha$ is surjective. Set $Q = J\cap T$. By the exactness of the previous sequence, $\Imag(\alpha)=\Tor_{2}^S(K,S/V)_{d+1}=\ker(\beta)$. Hence $\beta$ is the zero map and $\Tor^S_1(K,S/Q)_{d+1}=0$, as $\beta$ is surjective.

Therefore, if $\mathfrak{m}=(x_1,\dots,x_n)$,
$$
\Tor^S_1(K,S/Q)_{d+1}\cong\Tor^S_1(S/\mathfrak{m},S/Q)_{d+1}\cong\Big((Q\cap\mathfrak{m})\big/\mathfrak{m}Q\Big)_{d+1}=0,
$$
and $Q$ is generated only in one degree $d$. Since $G(I)=G(I)_d=[I_d]_t=[(J\cap T)_d]_t=[Q_d]_t=G(Q)_d=G(Q)$, $I$ and $Q=J\cap T$ are monomial ideals equigenerated in  degree $d$ and with the same minimal generating set. Hence, $I=Q$ and condition (b) holds. \\
The converse is also true, since if $I= J\cap T$, then $J \cap T$ is generated in degree $d$ and, using again the exactness of the sequence (\ref{exactsequencecomletelytspreadlexseg}), $\alpha$ will be surjective. 

Thanks to Claim 1, to get the equivalence (a)$\iff$(b), it suffices to prove that $\alpha$ is surjective if and only if (b) holds.

Set $\widetilde S=K[x_1,\dots,x_{n-(d-1)(t-1)}]$ and $\widetilde {\bf x}=x_1,\dots,x_{n-(d-1)(t-1)}$.

Let us identify the $K$--vector spaces $$\Tor^S_i(K,S/J),\ \ \Tor^S_i(K,S/T),\ \ \Tor^S_i(K,S/V)$$ with the Koszul homologies 
$$H_i(\widetilde{\bf x};\widetilde S/J^{\tau^{t-1}}),\ \ H_i(\widetilde{\bf x};\widetilde S/T^{\tau^{t-1}}),\ \ H_i(\widetilde{\bf x};\widetilde S/V^{\tau^{t-1}}),$$ respectively (see Lemma \ref{lemmakoszulhomtspread}, proof). In order to simplify the notations, we set 
$$\widetilde J=J^{\tau^{t-1}},\quad \widetilde T=T^{\tau^{t-1}},\quad \widetilde{V}=V^{\tau^{t-1}},$$
and
$$T_2(\widetilde J) = H_2(\widetilde{\bf x};\widetilde S/\widetilde J),\quad T_2(\widetilde T) =H_2(\widetilde{\bf x};\widetilde S/\widetilde T),\quad T_2(\widetilde V) =H_2(\widetilde{\bf x};\widetilde S/\widetilde V).$$

Hence, we must prove that
$$
\alpha:T_2(\widetilde J)\oplus T_2(\widetilde T)\rightarrow T_2(\widetilde V)
$$
is surjective if and only if condition (b) holds.

By Lemma \ref{AraHerHibiKoszulHomology}, the homology classes of the Koszul cycles (\emph{classical cycles})
	$$
	\omega/x_{\max(\omega)}e_i\wedge e_{\max(\omega)}
	$$
with $\omega\in G(\widetilde V)=G(I_{n,d,t}^{\tau^{t-1}})=G(I_{n-(d-1)(t-1),d,1})=M_{n-(d-1)(t-1),d,1}$, $i\notin\supp(\omega)$, $1\le i<\max(\omega)$, form a $K$--basis of $T_2(\widetilde V)$.
Analogously, homology classes of the same kind, with $\omega\in G(\widetilde{J})$, form a $K$--basis of $T_2(\widetilde J)$.\\
The morphism  $\alpha$ sends these last kind of homology classes to the corresponding homology classes in $T_2(\widetilde V)$. Hence, all elements $[\omega/x_{\max(\omega)}e_i\wedge e_{\max(\omega)}]$, with $\omega\in G(\widetilde{J})$, $1\le i<\max(\omega)$ and $i\not\in\supp(\omega)$ are in the image of $\alpha$.\\

Let $W$ be the $K$--vector subspace of $T_2(\widetilde V)$ generated by the following homology classes:
\begin{enumerate}
	\item[(I)] $[\omega/x_{\max(\omega)}e_i\wedge e_{\max(\omega)}]$, with $\omega\in\mathcal{L}_1^i(v^{\tau^{t-1}})=G(\widetilde{J})$ in the suitable polynomial ring $K[x_1,\dots,x_{n-(d-1)(t-1)}]=\widetilde{S}$, $1\le i<\max(\omega)$, $i\notin\supp(\omega)$;
	\item[(II)] $[\omega/x_{\max(\omega)}e_i\wedge e_{\max(\omega)}]$, with $\omega\in M_{n-(d-1)(t-1),d,1}$, $i_1<i<\max(\omega)$, $i\notin\supp(\omega)$.
\end{enumerate}

We are going to prove that $\Imag(\alpha)$ contains $W$. \\
Indeed, we have just observed that the elements described in (I) are in the image of $\alpha$. So, it remains to show that also the elements described in (II) are in $\Imag(\alpha)$.\\
We can observe that $G(\widetilde{J})=\mathcal{L}_1^i(v^{\tau^{t-1}})$ in the polynomial ring $\widetilde{S}$. Hence, if one considers $\omega\notin G(\widetilde{J})$, that is $\omega<_{\slex}v^{\tau^{t-1
}}$, then the differential of the corresponding classical Koszul cycles is
\[\partial_2(\omega/x_{\max(\omega)}e_i\wedge e_{\max(\omega)})= \omega e_i-x_i\omega/x_{\max(\omega)}e_{\max(\omega)},
\]
with $i\notin\supp(\omega)$ and $1\le i<\max(\omega)$.\smallskip

Suppose $i>i_1$. Since $\omega<_{\slex}v^{\tau^{t-1}}\le_{\slex}u^{\tau^{t-1}}$, $\min(\omega)\ge\min(u^{\tau^{t-1}})$ and $i>i_1=\min(u)$, we have $x_i(\omega/x_{\max(\omega)})\le_{\slex}u^{\tau^{t-1}}$. Hence, $x_i(\omega/x_{\max(\omega)})\in G(\widetilde{T})=\widetilde{F}$, where $\widetilde{F}=F^{\tau^{t-1}}$. Therefore, $x_i(\omega/x_{\max(\omega)})$, $\omega\in\widetilde{T}$, and
$$
\partial_2(\omega/x_{\max(\omega)}e_i\wedge e_{\max(\omega)})=0.
$$
Hence, $[\omega/x_{\max(\omega)}e_i\wedge e_{\max(\omega)}]\in T_2(\widetilde{T})$, for all monomials $\omega<_{\slex}v^{\tau^{t-1}}$, $i\notin\supp(\omega)$, $i_1<i<\max(\omega)$. And so, these homology classes are also in the image of $\alpha$. 

Each homology class $h=[\omega/x_{\max(\omega)}e_i\wedge e_{\max(\omega)}]$, with $\omega\in M_{n-(d-1)(t-1),d,1}$, $i_1<i<\max(\omega)$, $i\notin\supp(\omega)$ described in (II) belongs to $\Imag(\alpha)$. Indeed, if $\omega\ge_{\slex}v^{\tau^{t-1}}$, then the homology class $h$, we are keeping in consideration, is actually also described by (I), and so $h\in\Imag(\alpha)$. Otherwise, if $\omega<_{\slex}v^{\tau^{t-1}}$, then by the discussion above, $h\in\Imag(\alpha)$, too.\\\
Hence, $W\subseteq\Imag(\alpha)$ and thus we can consider the induced function
$$
\widetilde{\alpha}:T_2(\widetilde{T})\rightarrow T_2(\widetilde V)/W,
$$
and $\alpha$ is surjective if and only if $\widetilde{\alpha}$ is surjective. So, it suffices to prove that $\widetilde{\alpha}$ is surjective if and only if condition (b) holds.\\\\\
\textsc{Claim 2}. $\Imag(\widetilde{\alpha})$ is generated by the homology classes
\begin{equation}\label{basisImageAlphaKoszul}
[(zx_i/x_{i_1})/x_{\max(zx_i/x_{i_1})}e_{i_1}\wedge e_{\max(zx_i/x_{i_1})}],
\end{equation}
with $z\in\widetilde{T}=\mathcal{L}_1^f(u^{\tau^{t-1}})$ in the polynomial ring $\widetilde{S}$, $i\notin\supp(z)$ and $i_1<i\le\max(zx_i/x_{i_1})$.\\

Indeed, by (\ref{eq2:fintspread}), a basis of $T_2(\widetilde T)$ is given by the homology classes of the Koszul cycles $z/x_{\min(z)}e_{\min(z)}\wedge e_i$, with $z\in\mathcal{L}_1^f(u^{\tau^{t-1}})$ in $\widetilde S$, $\min(z)<i\le n-(d-1)(t-1)$, and $i\notin\supp(z)$. \\
We distinguish two cases depending on the value of $i$.\\
\textsc{Case 1.} Let $\max(z)< i\le n-(d-1)(t-1)$.\\ 
We have $\max(x_iz/x_{\min(z)})=i$ and, setting $\widetilde{z}=x_iz/x_{\min(z)}$, we can write
$$
z/x_{\min(z)}e_{\min(z)}\wedge e_i=\widetilde{z}/x_{\max(\widetilde{z})}e_{\min(z)}\wedge e_{\max(\widetilde{z})}.
$$
Since $z\le_{\slex}u^{\tau^{t-1}}$, we have that $\min(z)\ge\min(u^{\tau^{t-1}})=\min(u)=i_1$.

Now we distinguish two subcases.\\
\textsc{Subcase 1.} Let $\min(z)>i_1$. Then $\alpha([z/x_{\min(z)}e_{\min(z)}\wedge e_i])\in W$ and so
$$
\widetilde{\alpha}([z/x_{\min(z)}e_{\min(z)}\wedge e_i])=0.
$$
\textsc{Subcase 2.} Let $\min(z)=i_1$. We have that
\begin{align*}
\widetilde{\alpha}([z/x_{\min(z)}e_{\min(z)}\wedge e_i])&=[\widetilde{z}/x_{\max(\widetilde{z})}e_{\min(z)}\wedge e_{\max(\widetilde{z})}]\\
&=[(zx_i/x_{i_1})/x_{\max(zx_i/x_{i_1})}e_{i_1}\wedge e_{\max(zx_i/x_{i_1})}]
\end{align*}
is an element of $\Imag(\widetilde{\alpha})$ as described in \textsc{Claim 2.}\\\\
\textsc{Case 2.} Let $\min(z)<i<\max(z)$.\\
Let us  consider $z/(x_{\min(z)}x_{\max(z)})e_{\min(z)}\wedge e_i\wedge e_{\max(z)}$ $\in K_3(\widetilde{\bf x},\widetilde{S}/\widetilde{T})$. Then
\begin{align*}
\partial_3(z/(x_{\min(z)}x_{\max(z)})e_{\min(z)}&\wedge e_i\wedge e_{\max(z)})=z/x_{\max(z)}e_i\wedge e_{\max(z)}\\
&-x_iz/(x_{\min(z)}x_{\max(z)})e_{\min(z)}\wedge e_{\max(z)}\\
&+z/x_{\min(z)}e_{\min(z)}\wedge e_i.
\end{align*}
It follows that 
\begin{align*}
[z/x_{\min(z)}e_{\min(z)}\wedge e_i]&=[x_iz/(x_{\min(z)}x_{\max(z)})e_{\min(z)}\wedge e_{\max(z)}]\\&-[z/x_{\max(z)}e_i\wedge e_{\max(z)}].
\end{align*}
Note that the last summand is in $W$.

Recall that $\min(z)\ge i_1$.\\
If $\min(z)>i_1$, as before we have that $\widetilde{\alpha}([z/x_{\min(z)}e_{\min(z)}\wedge e_i])=0$.\\
If $\min(z)=i_1$, we observe that $\max(z)=\max(x_iz/x_{\min(z)})$ and thus
\begin{align*}
\widetilde{\alpha}([z/x_{\min(z)}e_{\min(z)}\wedge e_i])&=[x_iz/(x_{\min(z)}x_{\max(z)})e_{\min(z)}\wedge e_{\max(z)}]\\
&=[(x_iz/x_{\min(z)})/x_{\max(x_iz/x_{\min(z)})}e_{\min(z)}\wedge e_{\max(x_iz/x_{\min(z)})}]\\
&=[(zx_i/x_{i_1})/x_{\max(zx_i/x_{i_1})}e_{i_1}\wedge e_{\max(zx_i/x_{i_1})}]
\end{align*}
is an element of $\Imag(\widetilde{\alpha})$ as described in \textsc{Claim 2}.

Recalling that the generators of $W$ are those described in (I) and (II), we obtain that the canonical basis of $T_2(\widetilde V)/W$ is given by the homology classes of the Koszul cycles
\begin{equation}\label{basisCoDomainAlphaKoszul}
[\omega/x_{\max(\omega)}e_{i_1}\wedge e_{\max(\omega)}],
\end{equation}
with $\omega\in G(\widetilde V)=M_{n-(d-1)(t-1),d,1}$, $\omega<_{\slex}v^{\tau^{t-1}}$, $i_1\notin\supp(\omega)$.\\

Now, let us consider the following condition:
\begin{enumerate}
	\item[(b')] for every $\omega\in M_{n-(d-1)(t-1),d,1}$ with $\omega<_{\slex}v^{\tau^{t-1}}$ there is $i>i_1$ such that $i\in\supp(\omega)$ and $x_{i_1}(\omega/x_{i})\le_{\lex}u^{\tau^{t-1}}$.  
\end{enumerate}
We show that $\widetilde{\alpha}$ is surjective if and only if condition (b') holds.\smallskip

Indeed, assume (b') holds. \\ Let $h$ be an element of the type given in equation (\ref{basisCoDomainAlphaKoszul}). Then $$h=[\omega/x_{\max(\omega)}e_{i_1}\wedge e_{\max(\omega)}],$$ 
with $\omega\in G(\widetilde V)=M_{n-(d-1)(t-1),d,1}$, $\omega<_{\slex}v^{\tau^{t-1}}$, $i_1\notin\supp(\omega)$. By hypothesis, there is an integer $i>i_1$ such $i\in\supp(\omega)$ and $x_{i_1}(\omega/x_{i})\le_{\slex}u^{\tau^{t-1}}$. Hence, $\omega=zx_i/x_{i_1}$, with $z\in\mathcal{L}_1^f(u^{\tau^{t-1}})$. Therefore, since $\max(\omega)=\max(zx_i/x_{i_1})$, we have that
\begin{align*}
h&=[\omega/x_{\max(\omega)}e_{i_1}\wedge e_{\max(\omega)}]\\
&=[(zx_i/x_{i_1})/x_{\max(zx_i/x_{i_1})}e_{i_1}\wedge e_{\max(zx_i/x_{i_1})}]
\end{align*}
is an element given in (\ref{basisImageAlphaKoszul}) and so $\widetilde{\alpha}$ is surjective.\\
Suppose now that $\widetilde{\alpha}$ is surjective. Each element $h$ of the canonical basis of $T_2(\widetilde V)/W$ has a different multidegree. Hence, since $\widetilde{\alpha}$ is surjective, $h=\widetilde{\alpha}(h_1)$ for some homology class $h_1$ of a non--classical cycle of $T_2(\widetilde T)$. So, each element given in equation (\ref{basisCoDomainAlphaKoszul}) is equal to some element given in equation (\ref{basisImageAlphaKoszul}). Arguing as before, this implies that condition (b') holds.\\
Finally, the theorem is proved since the equivalence (b')$\iff$(b) follows from the next technical lemma.
\end{proof}

\begin{Lem}\label{lem:tec}
	Let $L=\mathcal{L}_t(u,v)$ and $I=(L)$ as in Theorem \ref{CharacterizationCompletelyLexSegTSpread}, then condition \textup{(b')} is equivalent to condition \textup{(b)}.
\end{Lem}
\begin{proof}
To prove the equivalence, we use Corollary \ref{cortaut-1lexseg} and the operator $$\sigma^{t-1}:M_{n-(d-1)(t-1),d,1}\rightarrow M_{n,d,t},$$ which is the inverse function of $\tau^{t-1}:M_{n,d,t}\rightarrow M_{n-(d-1)(t-1),d,1}$ \cite{EHQ}.\\\\
(b')$\Longrightarrow$(b). If $t=1$, then $\tau^{t-1}:M_{n,d,1}\rightarrow M_{n,d,1}$ is the identity map and there is nothing to prove. Therefore, we can suppose $t\ge2$.\\
Let $\omega =x_{k_1}x_{k_2}\cdots x_{k_d} \in M_{n,d,t}$ such that $\omega<_{\slex}v$. Then, $\omega^{\tau^{t-1}}\in M_{n-(d-1)(t-1),d,1}$ and $\omega^{\tau^{t-1}}=x_{k_1}x_{k_2-(t-1)}\cdots x_{k_{d}-(d-1)(t-1)}<_{\lex}v^{\tau^{t-1}}$ by Corollary \ref{cortaut-1lexseg}. By hypothesis, there exists an integer $i\in\supp(\omega^{\tau^{t-1}})$, $i=k_{\ell}-(\ell-1)(t-1)$, $\ell\in[d]$, such that $i>i_1$ and $x_{i_1}(\omega^{\tau^{t-1}}/x_{i})\le_{\lex}u^{\tau^{t-1}}$.
Since $\omega\le_{\slex}v\le_{\slex}u$, then $\min(\omega)\ge\min(u)=i_1$.
We have $k_1=\min(\omega)=\min(\omega^{\tau^{t-1}})>i_1$, otherwise $x_{i_1}^2$ divides $x_{i_1}(\omega^{\tau^{t-1}}/x_i)$ and consequently $x_{i_1}(\omega^{\tau^{t-1}}/x_i)>_{\lex}u^{\tau^{t-1}}$, a contradiction.\\
If $\ell=1$, then $x_{i_1}(\omega^{\tau^{t-1}}/x_i)=x_{i_1}x_{k_2-(t-1)}\cdots x_{k_d-(d-1)(t-1)}\le_{\lex}u^{\tau^{t-1}}$. Since $k_2-(t-1)-i_1> k_2-(t-1)-k_1\ge 1$, both monomials are squarefree and we can apply the operator $\sigma^{t-1}:M_{n-(d-1)(t-1),d,1}\rightarrow M_{n,d,t}$. We get $x_{i_1}x_{k_2}\cdots x_{k_d}=x_{i_1}(\omega/x_{k_1})\le_{\lex}u$, with $k_1>i_1$. Hence, condition (b) holds also in such a case.\\
If $1<\ell\le d$, then $i=k_{\ell}-(\ell-1)(d-1)$ and
$$
x_{i_1}(\omega^{\tau^{t-1}}/x_i)=x_{i_1}\Big(\prod_{r=1}^{\ell-1}x_{k_r-(r-1)(t-1)}\Big)\Big(\prod_{r=\ell+1}^{d}x_{k_r-(r-1)(t-1)}\Big)\le_{\lex}u^{\tau^{t-1}}.
$$
As before, both these monomials are squarefree. Applying $\sigma^{t-1}$, we get
$$
x_{i_1}\Big(\prod_{r=1}^{\ell-1}x_{k_r+(t-1)}\Big)\Big(\prod_{r=\ell+1}^{d}x_{k_r}\Big)\le_{\lex}u.
$$
Finally, since $\omega$ is $t$--spread, $k_2\ge k_1+t>k_1+(t-1)$ and we have
$$
x_{i_1}(\omega/x_{k_1})<_{\lex}x_{i_1}\Big(\prod_{r=1}^{\ell-1}x_{k_r+(t-1)}\Big)\Big(\prod_{r=\ell+1}^{d}x_{k_r}\Big)\le_{\lex}u,
$$
with $k_1>i_1$ as observed before.\\\\\
(b)$\Longrightarrow$(b'). Let $\omega=x_{k_1}x_{k_2}\cdots x_{k_d}\in M_{n-(d-1)(t-1),d,1}$, $\omega<_{\slex}v^{\tau^{t-1}}$. Then, $\omega^{\sigma^{t-1}}=x_{k_1}x_{k_2+(t-1)}$ $\cdots x_{k_d+(d-1)(t-1)}\in M_{n,d,t}$ and $\omega^{\sigma^{t-1}}<_{\slex}v$. From the hypothesis, for some $\ell\in[d]$, we have
$$
x_{i_1}\Big(\prod_{r=1}^{\ell-1}x_{k_r+(r-1)(t-1)}\Big)\Big(\prod_{r=\ell+1}^{d}x_{k_r+(r-1)(t-1)}\Big)\le_{\lex}u.
$$
If $\ell=1$, both monomials are $t$--spread since $k_2-i_1>k_2-k_1\ge t$, thus applying $\tau^{t-1}$, we obtain that $x_{i_1}(\omega/x_{k_1})\le_{\lex}u^{\tau^{t-1}}$.

Suppose $\ell>1$. Dividing by $x_{i_1}$.  We have that
$$
\Big(\prod_{r=1}^{\ell-1}x_{k_r+(r-1)(t-1)}\Big)\Big(\prod_{r=\ell+1}^{d}x_{k_r+(r-1)(t-1)}\Big)\le_{\lex}u/x_{i_1}=x_{i_2}\cdots x_{i_d}.
$$
Note that such monomials are $t$--spread. Applying $\tau^{t-1}:M_{n,d,t}\rightarrow M_{n-(d-1)(t-1),d,1}$, we get
$$
\Big(\prod_{r=1}^{\ell-1}x_{k_r}\Big)\Big(\prod_{r=\ell+1}^{d}x_{k_r+(t-1)}\Big)\le_{\lex}x_{i_2}x_{i_3-(t-1)}\cdots x_{i_d-(d-2)(t-1)}.
$$
Finally, multiplying by $x_{i_1}$ both sides, we have
\begin{align*}
x_{i_1}(\omega/x_{k_1})&\le_{\lex}x_{i_1}\Big(\prod_{r=1}^{\ell-1}x_{k_r}\Big)\Big(\prod_{r=\ell+1}^{d}x_{k_r+(t-1)}\Big)\\&\le_{\lex}x_{i_1}x_{i_2}x_{i_3-(t-1)}\cdots x_{i_d-(d-2)(t-1)}
\\&\le_{\lex}x_{i_1}x_{i_2-(t-1)}x_{i_3-2(t-1)}\cdots x_{i_d-(d-1)(t-1)}=u^{\tau^{t-1}},
\end{align*}
with $k_1>i_1$ as observed before and (b') holds.
\end{proof}

We collect some examples which illustrate condition (b) in Theorem \ref{CharacterizationCompletelyLexSegTSpread}.  
\begin{Expl}\rm\label{es:completely1}
Let $S=K[x_1,\dots,x_{12}]$ and $t=3$. Let $u=x_2x_6x_9$ and $v=x_3x_6x_9$. Set $L=\mathcal{L}_3(u,v)$.
We observe that $v\ne x_6x_9x_{12} = \min(M_{12, 3, 3})$. Thus $L$ is not a final $3$--spread lexsegment. Let $\omega=x_{j_1}x_{j_2}x_{j_3}<_{\slex}v$, then $j_1\ge 3,\ j_2\ge6$ and $j_3\ge10$. For such monomials, we have that, since $\min(\omega)>\min(u)$, $x_{\min(u)}\omega/x_{\min(\omega)}=x_2x_{j_2}x_{j_3}<_{\lex}u=x_2x_6x_9$, because $6\le j_2$ and $9<10\le j_3$. Hence, condition (b) in Theorem \ref{CharacterizationCompletelyLexSegTSpread} is satisfied and the ideal $I = (L)$ is a completely $3$--spread lexsegment ideal. Indeed, using \emph{Macaulay2}, we can check that $I=(\mathcal{L}_3^i(v))\cap(\mathcal{L}_3^f(u))= (\mathcal{L}_3(x_1x_4x_7,v))\cap(\mathcal{L}_3(u,x_6x_9x_{12}))$.
\end{Expl}
\begin{Expl}\rm
Let $S=K[x_1,\dots,x_{12}]$ and $d=3$, $t=2$. Let $u=x_2x_4x_6x_9$, $v=x_2x_4x_6x_{11}$ and set $L=\mathcal{L}_2(u,v)$.
We observe that $v\ne x_6x_8x_{10}x_{12}$, so $L$ is not final. In such a case condition (b) in Theorem \ref{CharacterizationCompletelyLexSegTSpread} is not verified. Indeed, we can consider the monomial $\omega=x_2x_4x_6x_{12}<_{\slex}v$. We have that $\min(u)=2$ and
\begin{align*}
x_{\min(u)}(\omega/x_4)&=x_2^2x_6x_{12}>_{\lex}u,\\
x_{\min(u)}(\omega/x_6)&=x_2^2x_4x_{12}>_{\lex}u,\\
x_{\min(u)}(\omega/x_{12})&=x_2^2x_4x_6>_{\lex}u.
\end{align*}
Hence, the ideal $I = (L)$ is not a completely $2$--spread lexsegment ideal. In fact, setting $J= (\mathcal{L}_2^i(v)) = (\mathcal{L}_2(x_1x_3x_5x_7,x_2x_4x_6x_{11}))$ and $T=(\mathcal{L}_2^f(u)) = (\mathcal{L}_2(x_2x_4x_6x_9,x_6x_8x_{10}x_{12}))$, let us consider the monomials $\omega_1=x_1x_4x_7x_{10}\in J$, $\omega_2=x_2x_4x_7x_{10}\in T$. It follows that the monomial $\textup{lcm}(\omega_1,\omega_2)=x_1x_2x_4x_7x_{10}$ belongs to the minimal system of monomial generators of $J\cap T$ but it does not belong to $G(I)$. Therefore, $I\ne J\cap T$.
\end{Expl}

\section{Completely $t$--spread lexsegment ideals with linear resolution}\label{sec3}

In this section, we want to characterize all completely $t$--spread lexsegment ideals with a linear resolution. The 1--spread case has also been considered in \cite{BS2008}.

Firstly, we need some comments and remarks.
Let $L=\mathcal{L}_t(u,v)$ be a $t$--spread lexsegment. Set $I=(\mathcal{L}_t(u,v))$. We may always assume
that $u\ne v$. Indeed, a principal ideal has a linear resolution. Furthermore, we may assume that $x_1$ divides $u$. Indeed, if $\min(u)>1$, then the sequence $x_1,\dots,x_{\min(u)-1}$ is regular on $S/I$, and the graded Betti numbers of $I$ are the same of $I'=I\cap K[x_{\min(u)},x_{\min(u)+1},\dots,x_n]$. Hence, $I$ has a linear resolution if and only if $I'$ has a linear resolution. Moreover, we may assume $\min(u)<\min(v)$. Indeed, if $\min(u)=\min(v)$, $I$ is isomorphic to the ideal $I''=(\mathcal{L}_t(u/x_{\min(u)},v/x_{\min(v)}))$. Hence, $I$ has a linear resolution if and only if $I''$ does. Finally, we assume that $I$ is not initial neither final, as in such cases, $I$ has a linear resolution.

Not all completely $t$--spread lexsegment ideals have a linear resolution. 
\begin{Expl}
	\rm 
	Let $I=(\mathcal{L}_0(x_1x_3^3,x_2^2x_3^2))=(x_1x_3^3,x_2^4,x_2^3x_3,x_2^2x_3^2)$ be a $0$--spread lexsegment ideal of the polynomial ring $\mathbb{Q}[x_1,x_2,x_3]$. Observe that $I$ is not final. By Lemma \ref{maxbettilemma}, for all $t\ge0$,
	\begin{align*}
	I^{\sigma^t}&=(\mathcal{L}_t(x_1x_{3+t}x_{3+2t}x_{3+3t},x_2x_{2+t}x_{3+2t}x_{3+3t}))\\
	&=(x_1x_{3+t}x_{3+2t}x_{3+3t},x_2x_{2+t}x_{2+2t}x_{2+3t},x_2x_{2+t}x_{2+2t}x_{3+3t},x_2x_{2+t}x_{3+2t}x_{3+3t})
	\end{align*}
	is a $t$--spread lexsegment of the polynomial ring $\mathbb{Q}[x_1,\dots,x_{3+3t}]$.\\
	For $t\ge 0$, since $\min(x_1x_{3+t}x_{3+2t}x_{3+3t})=1$ and $I^{\sigma^t}$ is not final, $I^{\sigma^t}$ is completely $t$--spread lexsegment if and only if for all $\omega\in M_{3+3t,4,t}$, $\omega<v$ it is $x_1(\omega/x_s)\le u$, for some $s\in\supp(\omega)$ (Theorem \ref{CharacterizationCompletelyLexSegTSpread} (b), \cite[Theorem 2.3 (b)]{DH}). Indeed, $\omega=x_2x_{3+t}x_{3+2t}x_{3+3t}$ or $\omega=x_3x_{3+t}x_{3+2t}x_{3+3t}$, and we may take $s=\min(\omega)$. Hence, for all $t\ge0$, $I^{\sigma^t}$ is a completely $t$--spread lexsegment ideal.

	Using \emph{Macaulay2}, one can verify that for $t=0,1,2, 3$, $I^{\sigma^t}$ has the following minimal resolution:
	
	\begin{equation}\label{eq:non}
	0\rightarrow R(-4)^2\oplus R(-5)\rightarrow R(-4)^4\rightarrow I^{\sigma^t}\rightarrow0,
	\end{equation}
	where $R=\mathbb{Q}[x_1,\dots,x_{3+3t}]$.
	Hence, in such cases $I^{\sigma^t}$ does not have a linear resolution. 
	
	Furthermore,  since $I$ is a monomial ideal of $\mathbb{Q}[x_1,x_2,x_3]$, set $J= I^{\sigma^3}$, one has that $J$ and $J^{\sigma^t}$, $t\ge 0$, have the same graded Betti numbers (\cite{EHQ}, Introduction). Indeed, for such ideals the application of the shifting operator $\sigma$ simply amounts to rename the variables. It follows, that (\ref{eq:non}) is a minimal graded resolution also for $I^{\sigma^t}$, with $t\ge 3$. Finally, $I^{\sigma^t}$, with $t\ge 0$ is a completely $t$--spread lexsegment that has no linear resolution.
\end{Expl}

As in the previous section, if $D$ is a homogeneous ideal of $S=K[x_1, \ldots, x_n]$, we denote by $T_i(D)$ the $K$--vector space $\Tor_i^S(K,S/D)\cong H_i({\bf x};S/D)$, where ${\bf x} = x_1, \ldots x_n$, for all $i\ge0$.
Moreover, in what follows, we write $J=(\mathcal{L}^i_t(v))$, $T=(\mathcal{L}_t^f(u))$, $V=I_{n,d,t}= J+T$ and $Q=J\cap T$.

\begin{Lem}\label{LemmaSurMapsTor}
	Let $t\ge0$ an integer and let $I=(\mathcal{L}_t(u,v))$ be a completely $t$--spread lexsegment ideal of degree $d$ of $S$. Then $I$ has a linear resolution if and only if the maps
	$$
	T_{i+1}(J)_j\oplus T_{i+1}(T)_j\xrightarrow{\ \ \ } T_{i+1}(V)_j
	$$
	are surjective, for all $i\ge1$ and all $j\ne i+d-1$.
\end{Lem}
\begin{proof}
	$I$ has a linear resolution if and only if $T_i(I)_j=0$, for all $i\ge1$ and all $j\ne i+d-1$. Let $i>0$ and $j\ne i+d-1$. Since $I$ is a completely $t$--spread lexsegment ideal, it is $I=J\cap T$. Hence, the long exact sequence (\ref{longexactsequencetortspread}) yields the long exact sequence of $K$--vector spaces
	$$
	\cdots \rightarrow T_{i+1}(J)_j\oplus T_{i+1}(T)_j\xrightarrow{\ \alpha\ } T_{i+1}(V)_j\xrightarrow{\ \beta\ } T_{i}(I)_j\xrightarrow{\ \gamma\ } T_{i}(J)_j\oplus T_{i}(T)_j \rightarrow \cdots,
	$$
	where the lower index $j$ denotes the $j$--th graded component of the corresponding vector space.
	
	Since $j\ne i+d-1$ and $J$ and $T$ have a linear resolution, then $T_i(J)_j=T_i(T)_j=0$, and one get the following exact sequence:
	$$
	\cdots \rightarrow T_{i+1}(J)_j\oplus T_{i+1}(T)_j\xrightarrow{\ \alpha\ } T_{i+1}(V)_j\xrightarrow{\ \beta\ } T_{i}(I)_j\xrightarrow{\ \ \ }0.
	$$
	By the exactness of the previous sequence, as $\beta$ is surjective, $T_i(I)_j=0$ if and only if $\beta$ is the zero map. Hence, if and only $\ker(\beta)=T_{i+1}(V)_j$. Again, by the exactness of the sequence, $\text{Im}(\alpha)=\ker(\beta)=T_{i+1}(V)_j$. This last condition is equivalent to the fact that $\alpha$ is surjective. Since $I$ is a completely $t$--spread lexsegment ideal the assertion follows from Theorem \ref{CharacterizationCompletelyLexSegTSpread}, proof.
	\end{proof}

For the proof of the next lemma, we need the following shifting operators:
\begin{align*}
\sigma^{t}:x_{i_1}x_{i_2}\cdots x_{i_d}\in M_{n,d,0}&\longmapsto\textstyle\prod_{j=1}^dx_{i_j+(j-1)t}\in M_{n+(d-1)t,d,t},\\
\tau^{t}:x_{i_1}x_{i_2}\cdots x_{i_d}\in M_{n,d,t}&\longmapsto\textstyle\prod_{j=1}^dx_{i_j-(j-1)t}\in M_{n-(d-1)t,d,0}.
\end{align*}
It is clear that if $I$ is a $0$--spread lexsegment ideal, then $I^{\sigma^t}$ is a $t$--spread lexsegment ideal (see, Corollary \ref{cortaut-1lexseg}).

\begin{Lem}\label{LemmaIandITau}
	Let $t\ge 0$ an integer and let $I=(\mathcal{L}_t(u,v))$ be a $t$--spread lexsegment ideal of $S=K[x_1,\dots,x_n]$ of degree $d$. 
	Assume $\min(u)=1$ and $\min(v)>1$. Then, $I$ is completely $t$--spread lexsegment ideal if and only if $I^{\tau^{t}}$ is a completely $0$--spread lexsegment ideal.
\end{Lem}
\begin{proof} Since an initial or a final $t$--spread lexsegment ideal is a completely $t$--spread lexsegment ideal, we may assume that $I$ is neither  initial nor final. 
We have $I^{\tau^t}=(\mathcal{L}_0(u^{\tau^t},v^{\tau^t}))$ in the polynomial ring $\widetilde{S}=K[x_1,\dots,x_{n-(d-1)t}]$. Set $\widetilde{n}=n-(d-1)t$. Setting $u^{\tau^t}=x_1^{a_1}x_2^{a_2}\cdots x_{\widetilde{n}}^{a_{\widetilde{n}}}$ and $v^{\tau^t}=x_1^{b_1}x_2^{b_2}\cdots x_{\widetilde{n}}^{b_{\widetilde{n}}}$, since $\min(u)=1$ and $\min(v)\ge 2$, we obtain $a_1\ne0$ and $b_1=0$. In particular, $a_1\ne b_1$. Hence, by \cite[Theorem 2.3 (b)]{DH}, $I^{\tau^t}$ is a completely 0--spread lexsegment ideal if and only if
	\begin{enumerate}
		\item[(i)] for every $\omega\in M_{\widetilde{n},d,0}$, $\omega<_{\lex}v^{\tau^t}$ there exists an integer $i>1$ such that $x_i$ divides $\omega$ and $x_1(\omega/x_i)\le_{\lex}u^{\tau^t}$.
	\end{enumerate}
	By Theorem \ref{CharacterizationCompletelyLexSegTSpread} (b), $I$ is a completely $t$--spread lexsegment ideal if and only if
	\begin{enumerate}
		\item[(ii)] for every $\omega\in M_{n,d,t}$, $\omega<_{\lex}v$ there exists an integer $s>i_1=1$ such that $x_s$ divides $\omega$ and $x_1(\omega/x_s)\le_{\lex}u$.
	\end{enumerate}
	Arguing as in Lemma \ref{lem:tec}, we obtain that (i) is equivalent to (ii), as desired.
\end{proof}

From Lemmas \ref{LemmaSurMapsTor} and \ref{LemmaIandITau}, we obtain a classification of all completely $t$--spread lexsegment ideals with a linear resolution.
\begin{Thm} \label{thm:linear}
	Let $I=(\mathcal{L}_t(u,v))$ be a completely $t$--spread lexsegment ideal of $S$ with $u=x_{i_1}x_{i_2}\cdots$ $ x_{i_d}$ and $v=x_{j_1}x_{j_2}\cdots x_{j_d}$. Assume $\min(u)=1$ and $\min(v)>\min(u)$. Then $I$ has a linear resolution if and only if one of the following conditions holds:
	\begin{enumerate}
		\item[\textup{(i)}] $i_2=1+t$;
		\item[\textup{(ii)}] $i_2>1+t$ and for the largest $\omega\in M_{n,d,t}$, $\omega<_{\lex}v$, we have $x_1(\omega/x_{\max(\omega)})\le_{\lex}x_1x_{i_2-t}x_{i_3-t}\cdots x_{i_d-t}$.
	\end{enumerate}
\end{Thm}
\begin{proof}
	By hypothesis $I=J\cap T$ with $J=(\mathcal{L}^i_t(v))$ and $T=(\mathcal{L}_t^f(u))$. Set $\widetilde{I}=I^{\tau^t}$, $\widetilde{J}=J^{\tau^t}$, $\widetilde{T}=T^{\tau^t}$ and $\widetilde{V}=V^{\tau^t}=I_{\widetilde{n},d,0}$, where $V=J+T=I_{n,d,t}$. It is well known that, for all $i$ and $j$, $T_i(J)_j\cong T_i(\widetilde{J})_j$, $T_i(T)_j\cong T_i(\widetilde{T})_j$, $T_i(V)_j\cong T_i(\widetilde{V})_j$. Hence, by Lemma \ref{LemmaSurMapsTor} and from these isomorphisms, $I$ has a linear resolution, if and only if the maps
	\begin{equation}\label{eq:CharCompletMapsSurj}
	T_i(\widetilde{J})_j\oplus T_i(\widetilde{T})_j\rightarrow T_i(\widetilde{V})_j
	\end{equation}
	are surjective, for all $i>0$ and all $j\ne i+d-1$. By Lemma \ref{LemmaIandITau}, $\widetilde{I}$ is a 0--spread completely lexsegment ideal, as $\min(u)=1$, $\min(v)\ge 2$ and $I$ is by hypothesis a completely $t$--spread lexsegment ideal. From (\ref{eq:CharCompletMapsSurj}), Lemma \ref{LemmaSurMapsTor} and the aforementioned isomorphisms, one has that $I$ has a linear resolution if and only if $\widetilde{I}$ has a linear resolution.
	
	Suppose $\widetilde{I}$ is a completely $0$--spread lexsegment ideal with a linear resolution. By \cite[Theorem 1.3]{ADH}, setting $u^{\tau^t}=x_1^{a_1}x_2^{a_2}\cdots x_{\widetilde{n}}^{a_{\widetilde{n}}}$ and $v^{\tau^t}=x_1^{b_1}x_2^{b_2}\cdots x_{\widetilde{n}}^{b_{\widetilde{n}}}$, we have that only one of the following conditions holds
	\begin{enumerate}
		\item[(a)] $u=x_1^px_2^{d-p}$ and $v=x_1^px_{\widetilde{n}}^{d-p}$, for some $1\le p\le d$;
		\item[(b)] $b_1<a_1-1$;
		\item[(c)] $b_1=a_1-1$ and for the largest $\widetilde{\omega}\in M_{\widetilde{n},d,0}$, $\widetilde{\omega}<_{\lex}v^{\tau^t}$, $x_1(\widetilde{\omega}/x_{\max(\widetilde{\omega})})\le_{\lex}u^{\tau^t}$.
	\end{enumerate}
	Condition (a) does not hold, as we would have $\min(u)=\min(v)=1$, against the assumption. To complete the proof it suffices to show that (b)$\iff$(i), and (c)$\iff$(ii).\medskip\\
	(b)$\Longrightarrow$(i). By hypothesis $\min(v)\ge 2$, so $b_1=0$. Therefore, $a_1-1>b_1$ is equivalent to the fact that $a_1\ge 2$, \emph{i.e.}, $x_1^2$ divides $u^{\tau^t}$. Applying $\sigma^t$, this last fact is equivalent to the fact that $x_1x_{1+t}$ divides $u$. Hence $i_2=1+t$, as desired. \\The converse quickly follows.\medskip\\
	(c)$\Longrightarrow$(ii). Indeed, $b_1=a_1-1$ implies $a_1=1$, so $i_2>1+t$. Moreover, let $\widetilde{\omega}$ be the largest monomial of $M_{\widetilde{n},d,0}$ such that $\widetilde{\omega}<_{\lex}v^{\tau^t}$. Then $\omega=\widetilde{\omega}^{\sigma^t}$ is the largest monomial of $M_{n,d,t}$ smaller than $v$. By condition (c), we have $x_1(\widetilde{\omega}/x_{\max(\widetilde{\omega})})\le_{\lex}u^{\tau^t}$. Since $\min(u)=1$, it follows that
	$$
	\widetilde{\omega}/x_{\max(\widetilde{\omega})}\le_{\lex}u^{\tau^t}/x_{\min(u)}=x_{i_2-t}x_{i_3-2t}\cdots x_{i_d-(d-1)t}.
	$$
	Thus, applying $\sigma^t$, we obtain $\omega/x_{\max(\omega)}\le_{\slex}x_{i_2-t}x_{i_3-t}\cdots x_{i_d-t}$,
	and condition (ii) follows.\medskip \\
	Conversely, assume (ii) holds. Then, $i_2>1+t$ implies $b_1=a_1-1$, as $b_1=0$. Moreover, $x_1(\omega/x_{\max(\omega)})\le_{\lex}x_1x_{i_2-t}x_{i_3-t}\cdots x_{i_d-t}$ implies $(\omega/x_{\max(\omega)})\le_{\lex}x_{i_2-t}x_{i_3-t}\cdots x_{i_d-t}$. Applying $\tau^{t}$, we get
	$$
	(\omega/x_{\max(\omega)})^{\tau^t}=\omega^{\tau^{t}}/x_{\max(\omega^{\tau^t})}\le_{\lex}x_{i_2-t}x_{i_3-2t}\cdots x_{i_d-(d-1)t}.
	$$
	Finally, multiplying by $x_1$, as $\widetilde{\omega}=\omega^{\tau^t}$, we have that (c) follows.
\end{proof}

We close this section with a formula for the graded Betti numbers.

\begin{Prop}
	Let $I=(\mathcal{L}_t(u,v))$ be a completely $t$--spread lexsegment ideal of degree $d$ of $S=K[x_1,\dots,x_n]$. Suppose $I$ has a linear resolution. Then, for all $i\ge0$,
	$$
	\beta_{i,i+d}(I)=\sum_{\omega\in\mathcal{L}_t^f(u)}\binom{n-\min(\omega)-(d-1)t}{i}-\sum_{\substack{\omega\in\mathcal{L}_t^f(v)\\ \omega\ne v}}\binom{\max(\omega)-(d-1)t-1}{i}.
	$$\normalfont
\end{Prop}
\begin{proof}
	Since $I$ has a linear resolution, $\beta_{i,i+j}(I)\neq 0$ 
	if and only if $j=d$. For all $i\ge 0$, a
	the long exact sequence (\ref{longexactsequencetortspread}) yields the exact sequence 
	$$
	\cdots \rightarrow T_{i+2}(V)_{\ell}\rightarrow T_{i+1}(Q)_\ell\rightarrow T_{i+1}(J)_\ell\oplus T_{i+1}(T)_\ell\rightarrow T_{i+1}(V)_\ell\rightarrow T_{i}(Q)_\ell\rightarrow \cdots, 
	$$
	for all $\ell \ge 0$, with $Q = J\cap T$. Let $\ell=i+d$. Since $\ell\ne i+d+1$, and $V$ has a linear resolution, $T_{i+2}(V)_\ell=0$. On the other hand, $I$ is a completely $t$--spread lexsegment ideal, then $Q=I$. Hence, since $\ell \neq i+d-1$, it follows that $\beta_{i,\ell}(S/Q)=\beta_{i,\ell}(S/I)=0$. Thus $T_{i}(Q)_\ell=0$ and we obtain the following short exact sequence
	$$
	0\rightarrow T_{i+1}(Q)_\ell\rightarrow T_{i+1}(J)_\ell\oplus T_{i+1}(T)_\ell\rightarrow T_{i+1}(V)_\ell\rightarrow 0.
	$$
	Therefore, from \cite{EHQ},
	\begin{align*}
	\beta_{i,i+d}(I)&=\beta_{i+1,i+d}(S/I)=\beta_{i+1,i+d}(S/J)+\beta_{i+1,i+d}(S/T)-\beta_{i+1,i+d}(S/V)\\
	&=\sum_{\omega\in G(J)}\binom{\max(\omega)-(d-1)t-1}{i}+\sum_{\omega\in G(T)}\binom{n-\min(\omega)-(d-1)t}{i}\\&-\sum_{\omega\in G(V)}\binom{\max(\omega)-(d-1)t-1}{i}.
	\end{align*}
	Finally, we get the desired formula, since $G(V)=G(J)\cup(\mathcal{L}_t^f(v)\setminus\{v\})$.
\end{proof}

\begin{Expl}\rm Let us consider the $3$--spread completely lexsegment ideal $I =(\mathcal{L}_t(u, v))$, with $u=x_1x_5x_8$, and $v=x_2x_5x_8$ in Example \ref{ex1}. 
$I$ has a linear resolution. Indeed, condition (ii) in Theorem \ref{thm:linear} is satisfied: $i_2=5>1+t=4$, $\omega=x_2x_5x_9$, $x_1(\omega/x_{\max(\omega)})=x_1x_2x_5\le_{\slex} x_{i_1}x_{i_2-t}x_{i_3-t}=x_1x_{5-3}x_{8-3}=x_1x_2x_5$.

Moreover, using \textit{Macaulay2}, the Betti tables of $S/I$, $S/J$, $S/T$ and $S/V$, with $J = (\mathcal{L}_t^i(v))$, $T= (\mathcal{L}_t^f(u))$, $V = J+T$,  are the following ones: 
	$$\begin{array}{cccccccc}
	\begin{matrix}
	&0&1&2&3&4\\\text{total:}&1&11&21&15&4\\\text{0:}&1&\text{.}&\text{.}&\text{.}&\text{.}\\\text{1:}&\text{.}&\text{.}&\text{.}&\text{.}&\text{.}\\\text{2:}&\text{.}&11&21&15&4\\\end{matrix}
	&&&&&&&
	\begin{matrix}
	&0&1&2&3&4&5\\\text{total:}&1&16&41&45&24&5\\\text{0:}&1&\text{.}&\text{.}&\text{.}&\text{.}&\text{.}\\\text{1:}&\text{.}&\text{.}&\text{.}&\text{.}&\text{.}&\text{.}\\\text{2:}&\text{.}&16&41&45&24&5\\\end{matrix} \\
	\textrm{\ \ \ Betti table of} \,\, S/I &&&&&&& \textrm{\ \ \ Betti table of} \,\, S/J\\ \\
	\begin{matrix}
	&0&1&2&3&4&5\\\text{total:}&1&30&85&96&50&10\\\text{0:}&1&\text{.}&\text{.}&\text{.}&\text{.}&\text{.}\\\text{1:}&\text{.}&\text{.}&\text{.}&\text{.}&\text{.}&\text{.}\\\text{2:}&\text{.}&30&85&96&50&10\\\end{matrix}
	&&&&&&&
	\begin{matrix}
	&0&1&2&3&4&5\\\text{total:}&1&35&105&126&70&15\\\text{0:}&1&\text{.}&\text{.}&\text{.}&\text{.}&\text{.}\\\text{1:}&\text{.}&\text{.}&\text{.}&\text{.}&\text{.}&\text{.}\\\text{2:}&\text{.}&35&105&126&70&15\\\end{matrix}\\
	\textrm{\ \ \ Betti table of} \,\, S/T &&&&&&& \textrm{\ \ \ Betti table of} \,\, S/V
	\end{array}$$
	Finally, one can verify that $$\beta_{i+1,i+d}(S/I)=\beta_{i+1,i+d}(S/J)+\beta_{i+1,i+d}(S/T)-\beta_{i+1,i+d}(S/V).$$	
\end{Expl}

\end{document}